 \documentclass[twoside]{article}

 \usepackage{iot}
 \usepackage{amssymb}
 \usepackage{amsmath}
 \newtheorem{theorem}{Theorem}[section]
 \newtheorem{lemma}[theorem]{Lemma}
 
 \newtheorem{rem}[theorem]{Remark}
 
 \newtheorem{prop}[theorem]{Proposition}
 \newtheorem{example}[theorem]{Example}

 \makeatletter
 \@addtoreset{equation}{section}
 \makeatother

 \makeatletter
 \def\section{\@startsection {section}{1}{\z@}{-1.5ex plus -.5ex
 minus -.2ex}{1ex plus .2ex}{\large\bf}}

 \newcommand{\F}{{\mathbb F}_2}
 \newcommand{\HA}[1]{{H_A(#1)}}
 \newcommand{\HB}[1]{{H_B(#1)}}
 \newcommand{\HC}[1]{{H_C(#1)}}
 \newcommand{\wt}[1]{{\widetilde{#1}}}
 \newcommand{\Hi}{{H_\infty}}
 \renewcommand{\b}[1]{{\overline{#1}}}
 \newcommand{\binomial}[2]{\genfrac{(}{)}{0pt}{}{#1}{#2}}

\def\PG(#1){\mathbb{PG}(#1)}
\def\AG(#1){\mathbb{AG}(#1)}

\ptitle{Complete caps in projective space which are disjoint from
a subspace of codimension two}
 \pauthor{David L. Wehlau}

\begin{document}
 \maketitle

 \begin{abstract}
  Working over the field of order 2 we consider those complete caps
which are disjoint from some codimension 2 subspace of projective
space.  We derive restrictive conditions which such a cap must
satisfy in order to be complete.  Using these conditions we
obtain explicit descriptions of complete caps which do not meet
every hyperplane in at least 5 points.  In particular, we
determine the set of cardinalities of all such complete caps in all
dimensions.
 \end{abstract}

\footnotetext{Corrected version, September 1, 2001}
\setcounter{footnote}{1}
\footnotetext{I
thank A. Davydov for pointing out a error
at the bottom of page 7,
 Sept.~15, 2003}

\section{Introduction}

  A subset of the projective space $\Sigma=\PG(n,2)$ is a {\it cap\/}
if no three points of $S$ are collinear.  A cap in $\Sigma$ is
said to be {\it complete\/} if it is not properly contained in any
other cap of $\Sigma$. A cap $S$ in $\Sigma$ is called {\it
large\/} if $|S| \geq 2^{n-1}+1$; otherwise $S$ is said to be
{\it small}. Much is known about the structure of large complete
caps.  In particular, it was shown in \cite{DT} that if $S
\subset \Sigma$ is a large complete cap then $|S| = 2^{n-1} +
2^i$ for some $i \in \{0,1,2,\ldots,n-3,n-1\}$.  In \cite{DT} it
is also shown that if $S$ is a cap in $\Sigma$ satisfying $|S|
\geq 2^{n-1} + 2$ then there exists a codimension 2 subspace $\Hi
\cong \PG(n-2,2)$ of $\Sigma$ such that $S \cap \Hi =
\emptyset$.  In \cite{BW} the same result was shown to hold for
all large caps. Conversely in \cite{BBW} it was shown that given
a fixed $r \geq 1$ then for all sufficiently large $n$ there
exists a cap $S$ in $\PG(n,2)$ such that $S$ meets every
codimension $r$ subspace $H \cong \PG(n-r,2)$ of $\PG(n,2)$.

 Here we will consider complete caps, $S \subset \Sigma=\PG(n,2)$ with the
property that some codimension 2 subspace $\Hi \cong \PG(n-2,2)$
of $\Sigma$ satisfies the condition $S \cap \Hi = \emptyset$. Fix
such an $\Hi$ and denote the three hyperplanes of $\Sigma$ which
contain $\Hi$ by $K_A, K_B$ and $K_C$. Furthermore we write $H_A
:= K_A \setminus \Hi$, $H_B := K_B \setminus \Hi$ and $H_C := K_C
\setminus \Hi$.  Write $A := H_A \cap S = K_A \cap S$, $B := H_B
\cap S = K_B \cap S$ and $C := H_C \cap S = K_C \cap S$.  Finally
define $A' = H_A \setminus A$, $B' = H_B \setminus B$, and $C' =
H_C \setminus C$.

  If $C = \emptyset$ then $S$ is complete if and only if
$S = \Sigma \setminus K_C$.  Henceforth we will assume that
none of the three sets $A$, $B$, and $C$ is empty.

  Consider the subspace
$\widetilde{F} \cong \PG(r,2)$ of $K_C$ generated by $C$.  Put $F
:= \widetilde{F} \cap \Hi$ and $\widehat{C} := \wt{F} \setminus
F$.  Then $F \cong \PG(r-1,2)$ and $\widehat{C} \cong \AG(r,2)$.
Consider the decomposition of $H_A$ and of $H_B$ into cosets of
$\widetilde{C}$. Each such coset $F'$ consists of $2^r$ points
such that $F' \sqcup F \cong \PG(r,2)$. There are $2^{n-r-1}$ such
cosets in each of $H_A$ and $H_B$. Each of these cosets is
isomorphic to $\wt{F} \setminus F \cong \AG(r,2)$. We will denote
the cosets in $H_A$ by $\HA{1}, \HA{2}, \dots, \HA{{2^{n-r-1}}}$.
Similarly we denote the cosets in $H_B$ by $\HB{1}, \HB{2},
\dots, \HB{{2^{n-r-1}}}$ ordered so that $\HB{i}=\widehat{C} +
\HA{i}$ for all $i$.  Write $A(i) := \HA{i} \cap A$, $A'(i) :=
\HA{i} \cap A'$, $B(i) := \HB{i} \cap B$ and $B'(i) := \HB{i}
\cap B'$.

  The cap $S$ is complete if and only if
\begin{eqnarray}
   A + B = C',\quad A + C = B',\quad B + C = A'; \\ \label{eq 1}
 (A \oplus A) \sqcup (B \oplus B) \sqcup (C \oplus C) = \Hi. \label{eq 2}
\end{eqnarray}

 In particular, if $S$ is complete we must have
\begin{equation} \label{coset equations}
A(i) + C = B'(i) {\rm \ and\ } B(i) + C = A'(i)
\end{equation}
for all $i=1,2,\ldots,2^{n-r-1}$.

  One way to satisfy Equation~\ref{coset equations} is to take
$A(i) = \HA{i}$ and $B(i) = \emptyset$.  Similarly we may take
$A(i) = \emptyset$ and $B(i) = \HB{i}$.  We will call these two
solutions the {\it trivial\/} solutions.  A second way to satisfy
this equation is to take $A(i) = \{a_1\}$ and $B(i) = \HB{i}
\setminus (a_1 + C)$.  Symmetrically we may take $B(i) = \{b_1\}$
and $A(i) = \HA{i} \setminus (b_1 + C)$.  If either $|A|=1$ or
$|B|=1$ then we call the solution a singleton solution.

  Once we have solutions to Equation~\ref{coset equations} we
still will need the extra conditions of $A + B = C'$ together
with Equation~\ref{eq 2}.  In Section~\ref{tangent hyp} we
explicitly describe all complete caps with $|C|=1$.  Using these
results we are able in Section~\ref{general case} to give
detailed sufficient conditions on a cap $S$ which satisfies
Equation~\ref{coset equations} for all $i$, for $S$ to be
complete.

In Sections~\ref{double plane}, \ref{triple threat} and
 \ref{Four Slice} we consider complete caps which meet at least one
hyperplane in less than 5 points. For such caps we are able to
find all solutions to Equation~\ref{coset equations} and thus to
explicitly describe such complete caps.

  In particular, we are able to construct many small complete caps.
Also note that the Plotkin doubling construction (described in
 Section~\ref{prelims}), as well as its generalization the Black/White
lift (described in \cite{Black/White lift}) can be used to build a great
many more small complete caps from the small complete caps constructed
here.

 \section{Preliminaries} \label{prelims}
Let $\F^{n+1}$ denote the $(n+1)$-dimensional vector space over
$\F$, the field of order 2.  The elements of $\PG(n,2)$ are the
one dimensional subspaces of $\F^{n+1}$.  Each such subspace may
be uniquely represented by the non-zero vector it contains. Fix a
basis $\{e_0, e_1, \ldots, e_n\}$ of $\F^{n+1}$ and consider a
point $x \in \PG(n,2)$.  For ease of notation, we write $x = s_1
s_2 \cdots s_r$ if $e_{s_1} + e_{s_2} + \ldots + e_{s_t}$ is the
non-zero vector in $x$ where $s_1, s_2, \ldots, s_r$ are distinct
elements of $\{0,1,\ldots,n\}$.  For example, $013$ denotes the
element of $\PG(n,2)$ comprised by the one-dimensional subspace
containing $e_0 + e_1 + e_3$.  From this point of view, the three
points of $\PG(n,2)$ represented by $x,y$ and $z$ are collinear if
and only if $x$, $y$ and $z$ lie in a plane in $\F^{n+1}$ which
occurs if and only if $x+y=z$ in $\F^{n+1}$. The usual inner
product on $\F^{n+1}$ induces an inner product on $\PG(n,2)$.  We
will write $(x)^\perp$ to denote the hyperplane of points in
$\PG(n,2)$ which are orthogonal to $x$ with respect to this inner
product. On occasion we will identify $\F^{n}$ with the affine
space $\AG(n,2)$.  In that setting we will use $\emptyset$ to
denote the zero vector in $\AG(n,2)=\F^{n}$.

  Let $X$ and $Y$ be two subsets of $\Sigma = \PG(n,2)$ and let
$z \in \PG(n,2)$.  We write $X \oplus Y$ to denote the set
$\{x+y \mid x \in X, y \in Y, x \neq y\} \subseteq \Sigma$.
If $X \cap Y = \emptyset$ we also write $X + Y = X \oplus Y$.
If $z \notin X$ we write $z + X := \{ z + x \mid x \in X\}
= \{z \} + X$.

  A subset $X$ of $\PG(n,2)$ is said to be {\it periodic\/} if there
exists a point $v \in \PG(n,2)$ such that $v + X = X$.  If such a
point $v$ exists it is called a {\it vertex\/} of $X$.  The
(possibly empty) set of all vertices of $X$ is denoted $V(X)$.
Note that if $X$ is periodic then necessarily $|X|$ is even.

  Let $X$ be a subset of $\Sigma = \PG(n,2)$.  Embed $\Sigma$ in
$\wt{\Sigma} \cong \PG(n+1,2)$ and let $v \in \wt{\Sigma}
\setminus \Sigma$.  We define the {\it Plotkin double\/} of $X$
(from the vertex $v$) by
$$\phi(X) := X \sqcup \{v+x \mid x\in X\}.$$  Then $\phi(X)$ is a
periodic subset of $\wt{\Sigma}$ with $V(\phi(X)) = \phi(V(X))
\sqcup \{v\}$.  It is straightforward to verify that $\phi(X)$ is a
cap in $\widetilde{\Sigma}$ if and only if $X$ is a cap in $\Sigma$;
also, $\phi(X)$ is complete if and only if $X$ is.

  A slightly more general version of the following lemma may be found
in \cite[Lemma~3.9]{BW}.  We include a proof here for the reader's
convenience.

\begin{lemma}\label{Sebastiens lemma}
  Let $S \subset \Sigma$ be complete.  Suppose there exists a hyperplane
$L$ of $\Sigma$ such that $S \cap L$ is periodic.  Then $S$ is
periodic. Furthermore, if $v$ is a vertex for $S \cap L$ then $v$
is also a vertex for $S$.
\end{lemma}
\begin{proof}
  Let $v$ be a vertex for $S \cap L$ and assume, by way of contradiction,
that $v$ is not a vertex for $S$.  Then there exists a point
$x_1$ of $S$ such that $x_2 := x_1 + v \notin S$.  Clearly $x_1
\notin L$ since $v \in L$ and thus $x_2 \notin L$.  Since $S$ is
complete, $x_2$ lies on a secant of $S$: say $x_2 = y_1 + z_1$
where $y_1 , z_1 \in S$.   Without loss of generality $y_1 \in L$
since every line of $\Sigma$ meets $L$. Therefore $y_2 := v + y_1
\in L \cap S$. But $y_2 = v + y_1 = v + x_2 + z_1 = x_1 + z_1$
and thus the line $\{y_2,x_1,z_1\}$ is fully contained in $S$, a
contradiction.
\end{proof}

\begin{lemma}\label{vertex subspace}
  Let $X$ be a subset of $\Sigma = \PG(n,2)$ and let $V(X)$ denote the
(possibly empty) set of vertices for $X$.  Then $V(X)$ is a projective
subspace of $\Sigma$.
\end{lemma}

\begin{proof}
  We must prove that $V(X)$ is closed under addition,
i.e., if $u,v$ are two distinct points of $V(X)$ then $u+v \in V(X)$.
Accordingly let $u,v$ be two distinct points of $V(X)$.
We must show that $w = v+u$ is also a vertex of $X$. Let $p \in
X$.  Since $u$ is a vertex, $u+p \in X$. Similarly, since $v \in
V(X)$ we have $w + p = v + (u+p) \in X$. Hence $w+X \subseteq X$ and
thus $w$ is a vertex of $X$.
\end{proof}

\section{Caps having a tangent hyperplane} \label{tangent hyp}

  In this section we characterize those complete caps $S$
which have a tangent hyperplane.  Many of the results of this
section are related to results found in \cite{BW}.

We suppose that $K_C$ is a hyperplane which is tangent to $S$,
i.e., that $C$ consists of a single point, $C = \{c_o\}$. Choose
a codimension 2 subspace $\Hi \cong \PG(n-2,2)$ contained in
$K_C$ and disjoint from $S$.  Clearly if $S$ is complete we must
have $B' = c_0 + A$ and $A'=c_0 + B$. Thus $S$ is determined by
$A$ (or by $B$). We wish to investigate conditions on the set $A$
which guarantee that $S$ is complete.

 Define $E = \Sigma \setminus (S \sqcup (S \oplus S))$ and $T = S
\sqcup E$.

\begin{lemma}
  $T$ is a complete cap in $\Sigma = \PG(n,2)$.
\end{lemma}
\begin{proof}
  Assume, by way of contradiction that $T$ contains a line, $\ell$
Clearly $\ell$ must contain at least two points of $E$, say
$e,e_1$. Clearly $E \subset K_C$ and thus $|\ell \cap K_C| \geq
2$ from which it follows that $\ell \subset K_C$.  Thus $\ell$
contains at least one point of $\Hi$, $e$ say.  Therefore, $e+A =
A'$ and $e+B = B'$. Now if $e_1 \in \Hi$ then $e_1 + A = A'$ and
$e_1 + B = B'$.  Conversely if $e_1 \in H_C$ then $e_1 + A = B'$
and $e_1 + B = A'$.  In both of these cases we have $e + e_1 + (A
\sqcup B) = (A \sqcup B)$.  This shows that $e + e_1 \in (A
\sqcup B) \oplus (A \sqcup B) \subset S \oplus S$. Therefore $e +
e_1 \notin T$.  This contradiction proves the lemma.
\end{proof}

\begin{prop}
  Let $A$ and $B$ be as defined above.
  Then $V(A) = V(B)$.
  In particular, $A$ is periodic if and only if $B$ is periodic.
\end{prop}

\begin{proof}
  Suppose $A$ is periodic, i.e., suppose there exists $v \in \Hi$ such
that $v + A = A$.    Since $v + H_A = H_A$ we have $v + A' = A'$.
Now $c_0 + B = A'$ and $c_0 + A = B'$. Therefore $v + B = v + c_o
+ A' = c_0 + v + A' = c_0 + A' = B$.
\end{proof}

  Next we give a sufficient condition on $A$ to guarantee that the
set $S = \{c_0\} \sqcup A \sqcup B$ is a {\it complete\/} cap.
The following theorem is similar to \cite[Theorem 4.1]{BW}.

\begin{theorem}\label{easy construction}
 Let $S$ be a cap in $\Sigma = \PG(n,2)$
with $|C|=1$ and $B' = c_0 + A$.  Suppose that $|A| \neq 2^{n-2}$.
If $A$ (or $B$) is not periodic then $S$ is a complete cap.
\end{theorem}

\begin{proof}
Clearly every point of $\Sigma \setminus L_C$ lies on a secant to
$S$ through $c_0$.  Let $\alpha = |A|$ and $\beta = |B|$.  Then
$\alpha + \beta = 2^{n-1}$. Since $\alpha \neq 2^{n-1}$ one of
$\alpha$ or $\beta$, say without loss of generality $\alpha$,
exceeds $2^{n-2}$.  Therefore $A \oplus A = \Hi$, i.e., every
point of $\Hi$ lies on a secant to $A$.

  Finally we prove that if $S$ is not complete (and
$|A| \neq 2^{n-2}$) then $A$ is periodic. Thus we suppose that
$S$ is not complete. By the foregoing this implies that there
exists $x \in H_C$ with $x \neq c_0$ such that $x \notin A+B$.
Then $x+A = B'$ and $c_0+A=B'$. Therefore $v+A=x+c_0+A=x+B'=A$.
Therefore $A$ is periodic with vertex $v$.
\end{proof}

\begin{rem} Note that the hypothesis that $|A| \neq 2^{n-2}$ is
required.  See Example~\ref{unique extender}.
\end{rem}

  Suppose that $S$ is a cap with $|C| = 1$ and $B' = c_0 + A$.
Then $|S|=2^{n-1}+1$.  Since $|S|$ is odd, $S$ cannot be periodic.
Thus by Lemma~\ref{Sebastiens lemma}, we have that both $A$ and
$B$ are not periodic.  Hence we have proved the following result.

\begin{theorem}\label{neat characterization}
  In $\Sigma = \PG(n,2)$, let $S$ be a cap which meets a hyperplane
$K_C$ of $\Sigma$ in a single point, $c_0$.  Choose a codimension
2 subspace $\Hi \cong \PG(n-2,2)$ contained in $K_C$ and disjoint
from $S$. Suppose $B' = c_0 + A$ and that $|A| \neq 2^{n-2}$.
Then $S$ is complete if and only if $A$ is not periodic if and
only if $B$ is not periodic.
\end{theorem}

  Next we consider in detail the case where $|A|=|B|=2^{n-2}$.

  Suppose first that $A \oplus A \neq \Hi$, i.e., suppose there exists
$p \in \Hi \setminus (A \oplus A)$.  Then $p + A = A'$.  Therefore
$A \oplus A = (p + A') \oplus (p + A') = A' \oplus A'
   = (c_0 + A') \oplus (c_0 + A') = B \oplus B$.
In particular, $p \notin B \oplus B$. Note that we have shown
that whenever $|A|=|B|=2^{n-2}$, we have $A \oplus A = B \oplus B
= A' \oplus A' = B' \oplus B'$.

\begin{lemma}
  If $|E \cap \Hi| \geq 2$ then $A, B, A'$ and $B'$ are each periodic.
\end{lemma}
\begin{proof}
  Let $e_1, e_2 \in E \cap \Hi$.  Then $e_1 + e_2 + A = e_1 + A' = A$ and
thus $A$ is periodic.  Similarly $e_1 + e_2$ is a vertex for $B,
A'$ and $B'$.
\end{proof}

  Suppose $A$ is not periodic.  By the preceding lemma, $|E| \leq 1$.
Hence either $S$ is complete or $T = S \sqcup \{e\}$.  Thus we
have proved the following theorem.

\begin{theorem}\label{best characterization}
 In $\Sigma = \PG(n,2)$, let $S$ be a cap which meets a hyperplane
$K_C$ of $\Sigma$ in a single point, $c_0$.  Choose a codimension
2 subspace $\Hi \cong \PG(n-2,2)$ contained in $K_C$ and disjoint
from $S$. Suppose $B' = c_0 + A$.   If $S$ is complete then $A$
and $B$ are not periodic. If $A$ (or $B$) is not periodic then
either
\begin{enumerate}
\item[\rm (i)] $S$ is complete  \quad or
\item[\rm (ii)] $|A|=|B|=2^{n-2}$ and $\Sigma = S \sqcup (S \oplus S)
\sqcup \{e\}$
     where $e \in \Hi$ and the complete cap $T = S \sqcup \{e\}$ is periodic
     with $V(T) = \{e + c_0\}$.
\end{enumerate}
\end{theorem}

\begin{rem} If $n \leq 4$ and $|A|=2^{n-2}$ then $A \oplus A \neq
\Hi$ since then $|A \oplus A| \leq \binomial{|A|}{2} < |\Hi|$.
\end{rem}

  The following four examples show that $A \oplus A$ may or may not equal
$\Hi$ and $A$ may or may not be periodic in all possible
combinations when $|A|=2^{n-2}$.

\begin{example}  
  In $\PG(9,2)$ take $A_0 = \{8\}$, $A_1 = \{08,18,28,38,48, 58,68,78\}$,
$A_5 = \{ijklm8 \mid 0 \leq i < j < k < l < m \leq 7\}$, $A_8 =
\{012345678\}$ and $A^\sharp = A_0 \sqcup A_1 \sqcup A_5 \sqcup
A_8$. With $K_A = (9)^\perp$ and $K_B = (8)^\perp$ we have
$A^\sharp \oplus A^\sharp = \Hi$ and $|A^\sharp| = 66$.  Now the
point $v := 01 \in \Hi$ lies on 16 secants to $A^\sharp$ and thus
$|A^\sharp \sqcup (v + A^\sharp)| = 116$. Hence we may easily
extend $A^\sharp$ to a set $A \subset (K_A \setminus \Hi)$ with
$|A|=128$ and $v+A = A$ by adding the 12 points of 6 new secant
lines of $v$ to $A^\sharp$.  Then $A \oplus A = \Hi$ and $A$ is
periodic.
\end{example}

\begin{example}  
  In $\PG(5,2)$ take $A = \{4,04,14,24,34,01234,0134,0234\}$ where
$K_A = (5)^\perp$ and $K_B = (4)^\perp$.  Then $A$ is not periodic
and $A \oplus A = \Hi$.  Thus $S = A \sqcup (45+A) \sqcup \{45\}$
is complete.
\end{example}

\begin{example}  
  In $\PG(5,2)$ take $A = \{4,04,14,014,24,024,34,034\}$ where
$K_A = (5)^\perp$ and $K_B = (4)^\perp$.  Then $A$ is periodic
(with vertex 0) and $A \oplus A \neq \Hi$. Here if $c_0=45$ then
$E = \{045,123,0123\}$.
\end{example}

\begin{example} \label{unique extender} 
  In $\PG(4,2)$ take $A = \{3,03,13,23\}$ where
$K_A = (4)^\perp$ and $K_B = (3)^\perp$.  Then $A$ is not periodic
and $A \oplus A \neq \Hi$.  Thus $\Hi \setminus (A \oplus A) =
\{e=012\}$ and $S = A \sqcup (c_0 + A) \sqcup \{c_0,e\}$ is a
complete cap.
\end{example}


\section{General Case} \label{general case}
  Next we suppose that the sets $A$ and $B$ are such that
Equation~\ref{coset equations} is satisfied for all $i$. Under
this hypothesis we want to investigate conditions on $A$ and $B$
which are sufficient to guarantee that the $S$ is a {\it
complete\/} cap.

  We suppose that
$\{1,2,\dots,t\} = \{i \mid \HA{i} \subset S\}$ for some $t$ with
$0 \leq t \leq 2^{n-r-1}$. Thus $\{1,2,\dots,t\} = \{i \mid
\HB{i} \cap S = \emptyset\}$. Further suppose that
$\{t+1,t+2,\dots,t+u\} = \{i \mid \HA{i} \cap S = \emptyset\}$
for some $u$ with $0 \leq u \leq 2^{n-r-1}-t$. Then
$\{t+1,t+2,\dots,t+u\} = \{i \mid \HB{i} \subset S\}$.  Thus the
pair of cosets corresponding to $i=1,2,\ldots,t+u$ are precisely
those having a trivial solution to Equation~\ref{coset equations}.

  Let $\b{\Sigma} \cong \PG(n-r,2)$ denote the quotient geometry
of $\Sigma$ by $F$.  This is just the projective geometry
associated to the quotient of the vector space associated to
$\Sigma=\PG(n,2)$ by the vector space associated to $F$.

Let $\b{a_i}$ and $\b{b_i}$ denote the points of $\b{\Sigma}$
corresponding to the cosets $\HA{i}$ and $\HB{i}$ respectively.
Write $\b{C}$ for the point of $\b{\Sigma}$ corresponding to
$\wt{F}$.  Similarly we define $\b{H_A}$, $\b{H_B}$, etc.\ for
the subsets of $\b{\Sigma}$ corresponding to $H_A$, $H_B$, etc.\
respectively.

Corresponding to the full, empty and non-empty cosets of $A$ and
$B$ we define the following subsets of $\b{\Sigma}$. Put $\b{A_f}
= \{\b{a_1},\b{a_2},\dots,\b{a_t}\}$, $\b{B_f} =
\{\b{b_{t+1}},\b{b_{t+2}},\dots,\b{b_{t+u}}\}$, $\b{A_e} =
\{\b{a_{t+1}},\b{a_{t+2}},\dots,\b{a_{t+u}}\}$, $\b{B_e} =
\{\b{b_1},\b{b_2},\dots,\b{b_t}\}$, $\b{A_{ne}} = \b{H_A}
\setminus \b{A_e}$ and $\b{B_{ne}} = \b{H_B} \setminus \b{B_e}$.
Recall that $\widehat{C} = \wt{F} \setminus F$.

\begin{theorem} \label{sufficient conditions}
  Suppose that $A$ and $B$ are such that Equation~\ref{coset equations}
is satisfied for all $i =1 ,2 \ldots, 2^{n-r-1}$.  Suppose $0 <
t+u < 2^{n-r-1}$, that $t \neq u$, and that either $\b{A_f}$ or
$\b{B_f}$ is non-periodic.  Further suppose that $\widehat{C}
\setminus C \subseteq \cup_{i=1}^{2^{n-r-1}} \left(A(i) +
B(i)\right)$. Then $S$ is a complete cap in $\Sigma$.
\end{theorem}

\begin{proof}
  By construction, $S$ is a cap in $\Sigma$.  Also Equation~\ref{coset equations}
guarantees that every point of $A'$ and every point of $B'$ lies
on at least one secant of $S$ through a point of $C$. Thus we
need to show that every point of $\Hi$ and every point of $H_C
\setminus C$ lies on a secant of $S$.

  Without loss of generality suppose that $\b{A_f}$ is
non-periodic.
  Consider the cap $\b{S_1}$ in $\b{\Sigma}$ given by
$\b{S_1} = \b{A_f} \sqcup \b{B_{ne}} \sqcup \{\b{C}\}$.  Applying
Theorem~\ref{best characterization} we see that for every point
$\b{c_k}$ of $\b{H_C}$ different from $\b{C}$ there exists $i \leq
t$ and $j \geq t+1$ such that $\b{a_i} + \b{b_j} = \b{c_k}$. Since
the coset $\HA{i}$ is entirely contained in $S$, and $\HB{j} \cap
S \neq \emptyset$, this means that every point of the coset
$\HC{k} = \HA{i} + \HB{j}$ lies on a secant to $S$.

By assumption, every point of $\widehat{C} \setminus C$ is
contained in $A + B$.

  Hence it only remains to prove that every point of $\Hi$ lies on
at least one secant to $S$.  Note that every point of $F$ lies in
$\HA{1} \oplus \HA{1}$ (or in $\HB{1} \oplus \HB{1}$ if $t=0$).
Thus we consider a point $z \in \Hi \setminus F$ and write $\b{z}$
for the point of $\b{\Sigma}$ corresponding to the coset of $F$
generated by $z$.  As above, if $\b{z} \in \b{A_f} + \b{A_{ne}}$
or $\b{z} \in \b{B_f} + \b{B_{ne}}$ then $z$ lies on a secant line
to $S$. Thus we assume, by way of contradiction, that $\b{z} +
\b{A_f} \subset \b{A_e}$ and $\b{z} + \b{B_f} \subset \b{B_e}$.
The first inclusion implies that $t \leq u$ while the second
implies that $u \leq t$.  Thus $t=u$ contradicting our hypothesis.
\end{proof}

\begin{rem}  \label{not all trivial}
  Note that if there exists at least one value of $i$ for which
a singleton solution occurs then $\widehat{C} \setminus C
\subseteq \cup_{i=1}^{2^{n-r-1}} \left(A(i) + B(i)\right)$.
Conversely if $t+u=2^{n-r-1}$ then $A(i) + B(i) = \emptyset$ for
all $i$ and thus $\widehat{C} \setminus C \not\subseteq
\cup_{i=1}^{2^{n-r-1}} \left(A(i) + B(i)\right)$. Hence if $t+u
=2^{n-r-1}$ then $S$ is not complete.  Furthermore, in order that
$1 \leq t+u \leq 2^{n-r-1}-1$ we must have $r \leq n-2$.
\end{rem}

\begin{rem} \label{t or u odd}
Note that, if $t$ or $u$ is odd, then $\b{A_f}$ or $\b{B_f}$
respectively is non-periodic.  Thus for a fixed value of $t+u$
with $1 \leq t+u \leq 2^{n-r-1}-1$ we may take $t=1$, for example,
to arrange that $\b{A_f}$ is non-periodic.  Furthermore, if $t+u
\neq 2$ we may simultaneously arrange that $t \neq u$.
\end{rem}

\begin{rem} There do exist complete caps for which $t+u = 0$
(and thus $\b{A_f} = \b{B_f} = \emptyset$).  See
Example~\ref{partition cap example}.  The condition that $t \neq
u$ is also sufficient but not necessary.
\end{rem}

\section{A family of examples} \label{family}

  In this section we consider the case where $|C| = 2^r - 1$,
i.e., where $\widehat{C} \setminus C$ is a single point $\{c_0\}$.
If $r=0$ then $C$ is empty contrary to our earlier assumption.  If
$r=1$ then $|C|=1$ and we are in the case considered in
Section~\ref{tangent hyp}.  Thus we will assume that $r \geq 2$.
Accordingly, we must have $n \geq 3$.

  Using the cap property of $S$ we see that if $|S \cap \HA{i}| \geq 2$
then $S \cap \HB{i} = \emptyset$.  Furthermore if $S \cap \HB{i}$
is a single point, $\{\alpha_i\}$, then using the completeness of
$S$ we must have $S \cap \HB{i} = \{\beta_i := \alpha_i + c_0\}$.
Finally if $S \cap \HB{i} = \emptyset$ then the completeness of
$S$ implies that $\HB{i} \subset S$.

Thus the only solutions to Equation~\ref{coset equations} are the
trivial solutions and the singleton solutions. Let $s$ denote the
number of these pairs of cosets which have the singleton solution.

Then $|S| = (2^r-1) + 2s + 2^r(2^{n-r-1}-s) = 2^{n-1} + 2^r -1 -
(2^r-2)s$.  By Remark~\ref{not all trivial}, $s$ cannot equal 0.
Furthermore by Theorem~\ref{sufficient conditions} (and
Remark~\ref{t or u odd}) there do exist complete caps of this form
for all $s = 1,2,\ldots,2^{n-r-1}-3$ and $s=2^{n-r-1}-1$. Thus we
find complete caps of this form of all cardinalities: $2^{n-r} +
k(2^r-2) + 1$ for $k=2,4,5,\ldots,2^{n-r-1}$ where $2 \leq r
\leq n-2$.

  Next we consider the case $s=2^{n-r-1}-2$.  Then $u+v=2$ and
either $u=v$ or else both $\b{A_f}$ and $\b{B_f}$ are periodic.
Reordering the cosets we may suppose $|H_A(i)|=|H_B(i)|=1$ for
$i=3,4,\ldots,2^{n-r-1}$. Write $\{\alpha_i\} = H_A(i)$ and
$\{\beta_i\} = H_B(i)$ for $i=3,4,\ldots,2^{n-r-1}$. We consider
first the case that $t=u=1$. Examining the proof of
Theorem~\ref{sufficient conditions} we see that every point of
$\PG(n,2)$ lies on a secant to $S$ except possibly the $2^r$
points in the coset $H_A(1) + H_A(2) = H_B(1) + H_B(2) =
\{z_1,z_2,\ldots,z_{2^r}\} \subset \Hi$.  These points $z_k$ can
only lie on secants of the form $\alpha_i + \alpha_j = z_k$ or
$\beta_i + \beta_j = z_k$.  Since $\beta_i + \beta_j = \alpha_i +
\alpha_j$ we may restrict our attention to the $\alpha_i$.  Thus
we see that $S$ is complete if and only if for every
$k=1,2,\ldots,2^r$ there exists $3 \leq i(k) < j(k) \leq
2^{n-r-1}$ such that $\alpha_{i(k)} + \alpha_{j(k)} = z_k$. Since
the $\alpha_i$ all lie in different cosets, it is clear that
$\alpha_{i(1)}, \alpha_{j(1)},\ldots, \alpha_{i(2^r)},
\alpha_{j(2^r)}$ must be distinct.  Thus we must have $2(2^r) \leq
2^{n-r-1}-2$ or equivalently $2r+3 \leq n$.

The case $(t,u)=(2,0)$ (or $(0,2)$) is entirely similar to the
above case.  The only difference being that the $2^r$ points which
may not lie on secants to $S$ form a coset in $H_C$ instead of in
$\Hi$.  In summary we see that there exists a complete cap $S$
with $|\widehat{C} \setminus C|=1 $ and with
$|S|=2^{n-r}+2^{r+1}+2^r-5$ (and $s=2^{n-r-1}-2$) if and only if
$n \geq 2r+3$.

\bigskip
  The only remaining possibility is that $s=2^{n-r-1}$, i.e.,
$t=u=0$.  Thus we now consider caps with $\widehat{C}\setminus C
= \{c_0\}$ and with $|\HA{i}| = |\HB{i}| = 1$ for all $i =
1,2,\ldots,2^{n-r-1}$.

\begin{lemma}
  Let $2 \leq r \leq n-1$ and suppose that $\widehat{C} \setminus C$
is a single point $\{c_0\}$ and that $|A(i)|=|B(i)|=1$ for all
$i=1,2,\ldots,2^{n-r-1}$.  Then $S$ is a complete cap if and only
if $A \oplus A = \Hi \setminus F$.
\end{lemma}

\begin{proof}
  Note that $A \oplus A = B \oplus B$ since
$\alpha_i + \alpha_j = (\alpha_i + c_0) + (\alpha_j + c_0) =
\beta_i + \beta_j$.  Furthermore $A \oplus A \subseteq \Hi
\setminus F$ and $C \oplus C = F$.  Thus $\Hi \subset S \oplus S$
if and only if $A \oplus A = \Hi \setminus F$.

  Suppose now that $A \oplus A = \Hi \setminus F$.  By construction
$A + C = B'$ and $B+C = A'$ thus $S$ is complete if and only if
$A+B = C'$.  Now $A + B = A + (c_0 + A) = \{c_0\} \sqcup \left(c_0
+ (A \oplus A)\right) = \{c_0\} \sqcup \left(c_0 + (\Hi \setminus
F)\right) = C'$. Hence we conclude that $S$ is complete if and
only if $A \oplus A = \Hi \setminus F$.
\end{proof}

 We will now attempt to determine those sets
$A = \{\alpha_1,\alpha_2,\ldots,\alpha_{2^{n-r-1}}\}$ for which
$A \oplus A = \Hi \setminus F$.  The remainder of this section
will be devoted to answering this question.  We will reduce this
question to an equivalent geometric question in $\b{K_A}$ where,
as above, $\b{K_A}$ denotes a hyperplane in $\b{\Sigma}$, the
quotient geometry of $\Sigma$ with respect to $F$.

  To do this we will fix some subspace
$L \cong \PG(n-r-1,2) \subset K_A$ such that $L \cap F =
\emptyset$. Choose an identification of one of the cosets
$\HA{i}$, say $\HA{i}$, with $\AG(r,2)$.  Choose this
identification so that the unique point $a_0 := \in \HA{1} \cap
L$ is the point identified with the zero vector in $\AG(r,2)$.
Note that with this identification if $u, v \in \AG(r,2)$
correspond to $x, y \in \HA{1}$ respectively then $u+v \in
\AG(r,2)$ corresponds to $x+y+a_0 \in \HA{1}$.

Now for each point $a$ of $K_A \setminus F$ there exists a unique
point $x$ of $\AG(r,2)$ such that $x + a \in L$. For each $x \in
\AG(r,2)$ define $(K_A)_x := \{ a \in (K_A \setminus F) \mid x+a
\in L\}$.  This partitioning of $K_A \setminus F$ induces a
partitioning of $A$: $A = \sqcup_{x \in \AG(r,2)} A_x$ where $A_x
:= A \cap (K_A)_x$. Similarly we have the partitioning of $\Hi$:
$\Hi = \sqcup_{x \in \AG(r,2)} (\Hi)_x$ where $(\Hi)_x := \Hi
\cap (K_A)_x \cong \PG(n-r-2)$ for all $x \in \AG(r,2)$.

\begin{example}\label{partition example}
 Take $n=7$ and $r=2$.  Put $K_A = (7)^\perp$,
$K_B = (6)^\perp$ and $K_C = (67)^\perp$.  Then $\Hi \cong
\PG(5,2)$ is the projective subspace generated by
$\{0,1,2,3,4,5\}$.    Take $C = \{067,167,0167\}$ so that $c_0
=67$.  Then $\widehat{C} = \{067,167,0167,67\}$, $\wt{F} \cong
\PG(2,2)$ is generated by $\{0,1,67\}$ and $F$ is the line
$\{0,1,01\}$. Choose $L \cong \PG(4,2)$ to be the subspace of
$K_A$ generated by $\{2,3,4,5,\allowbreak6\}$.
Then $L \cap F = \emptyset$
as required. We take $\HA{1}$ to be the coset containing $6$
say.  Thus $\HA{1} = \{6,06,16,016\}$ and $a_0 = 6 \in \HA{1}
\cap L$.  Then for example if we choose our identification so
that $u= 1 \in \AG(2,2)$ corresponds to $16 \in \HA{1}$ we have
$(\Hi)_u = (\Hi)_{16} =
\{1,12,13,14,15,123,124,125,134,135,145,1234,\allowbreak
1235,1245,1345,12345\}$.
\end{example}

  Now the useful property of this partition is that
for $u,w \in \AG(r,2)$ we have $A_u \oplus A_v \subseteq (\Hi)_{u+v}$.
Thus $A \oplus A = \Hi \setminus F$ if and only if
for every $w \in \AG(r,2)$ we have
$\cup_{u \in \AG(r,2)} \left( A_u \oplus A_{u+w} \right) =
  (\Hi)_w \setminus F$.
Note that the partition of $A$ yields a partition of $\b{A} =
\b{H_A}$: $$\b{H_A} = \sqcup_{x \in \AG(r,2)} \b{A_x}.$$
Furthermore $A \oplus A = \Hi \setminus F$ if and only if for
every $w \in \AG(r,2)$ we have
\begin{equation} \label{partition condition}
\cup_{u \in \AG(r,2)} \left( \b{A_u} \oplus \b{A_{u+w}} \right) =
  \b{(\Hi)}.
\end{equation}

Thus we have proved the following proposition:
\begin{prop} \label{partition criterion}
  Let $C \cong \AG(r,2) \setminus \{c_0\}$ and suppose $|A(i)| =
  |B(i)| = 1$ for all $i = 1,2,\ldots,2^{n-r-1}$.  Then $S$ is a
  complete cap if and only if the partition of $\AG(n-r,2) =
  \b{H_A}$ induced by $A$ satisfies Equation~\ref{partition condition} for all
  $w \in \AG(r,2)$.
\end{prop}

\begin{example} \label{partition cap example}
  We continue with the notation of the previous example.
Thus $\b{H_A} \cong \AG(4,2) \subset \b{K_A} \cong \PG(4,2)$.
Hence to obtain a complete cap $S$ we require a partition of
$\AG(4,2)$ into 4 subsets,
$\b{A_0},\b{A_1},\b{A_{01}},\b{A_{\emptyset}}$ satisfying
Equation~\ref{partition condition} for all $w \in \AG(2,2) =
\{0,1,01,\emptyset\}$. It is not too difficult to find such
partitions.  For example,
          $\b{A_{\emptyset}} := \emptyset$,
          $\b{A_0} := \{6,246,346,256,2346,2456\}$
          $\b{A_1} := \{236,356,456,2356,3456\}$,
and    $\b{A_{01}} := \{26,36,46,56,23456\}$. This partition
corresponds to the complete cap $S = A \sqcup B \sqcup C$ where
which $A = \{1236,1356,1456,12356,13456,06,0246,0346,0256,02346,
02456,0126,0136,\allowbreak 0146,\allowbreak 0156,0123456\}$ and $B = c_0 + A
= 01 + A$.
\end{example}

\begin{lemma} 
 Let $\AG(k,2) \subset \Lambda = \PG(k,2)$ denote
the complement of a hyperplane $H$ of $\Lambda$.  Suppose there
exists a partition of $\AG(k,2)$ indexed by $\AG(r,2): \AG(k,2) =
\sqcup_{w \in \AG(r,2)} X_w$ such that
$$
 \hskip 7em
\cup_{u \in \AG(r,2)} \left( \b{X_u} \oplus \b{X_{u+w}}\right)=
H  \hskip 7em(\ref{partition condition})
$$
 for all $w \in \AG(r,2)$.  Then  $r \leq k-2$.
\end{lemma}

\begin{proof}
  The partition of $\AG(k,2)$ induces a partition of the secant
lines of $\AG(k,2)$ as follows.  If $a_1 \in X_u$ and $a_2 \in
X_v$ then we say the secant line through $a_1$ and $a_2$ is of
type $u+v \in \AG(r,2)$.  Equation~\ref{partition condition} is
the condition that every point of $H$ lies on at least one secant
line of type $w$ for all $w \in \AG(r,2)$. Since $|H|=2^k-1$ and
since there are $\binomial{2^k}{2} = 2^{k-1}(2^k-1)$ secant lines
we must have $2^{k-1}(2^k-1) \geq 2^r(2^k-1)$. Thus $k-1 \geq r$.
Furthermore, if $r = k-1$ then every point of $H$ must lie on
exactly one secant line of each type $w \in \AG(r,2)$.

  Suppose then that $r=k-1$ and assume by way of contradiction that
a partition of $\AG(k,2)$ satisfying Equation~\ref{partition
condition} exists.  Define $x_w := |X_w|$ for $w \in \AG(k,2)$.
Then we have $\sum_{u \in \AG(r,2)} x_w x_u = 2^k-1$ for all $w$
different from zero in $\AG(r,2)$.  Also $\sum_{u \in \AG(r,2)}
\binomial{x_u}{2} = 2^k-1$.  Since $\sum_{u \in \AG(r,2)} x_u =
2^r-1$, this gives $\sum_{u \in \AG(r,2)} x_u^2 = 2^k -2$.  Thus
\samepage
\begin{eqnarray*}
2^{2k} &=& (\sum_{u \in \AG(r,2)} x_u)^2 = \sum_{u
\neq v} x_u x_v + \sum_{u \in \AG(r,2)} x_u^2\\
       &=& 2^r(2^k-1) +
2^k-2 = 2^{k-1}(2^k-1) + (2^k-1) - 1 = (2^k-1)(2^k).
\end{eqnarray*}
This contradiction shows that $r \leq k-2$.
\end{proof}

\begin{rem}
Clearly, if a partition of $\AG(k,2)$ exists which satisfies
Equation~\ref{partition condition} for some value $r_0$ of $r$
then such a partition exists for all values of $r$ less than
$r_0$.
\end{rem}

\section{$|C|=2$}\label{double plane}

\begin{lemma}
  Let $S$ be a complete cap in $\Sigma = \PG(n,2)$ satisfying
$|C|=2$.  Then $S = \phi(S')$ where $S'=S \cap \Sigma'$ is a
complete cap in a hyperplane $\Sigma' \cong \PG(n-1,2)$.  In
particular, $|S| = 2^{n-1}+2$.
\end{lemma}

\begin{proof}
  Write $C=S \cap K_C=\{c_1,c_2\}$.  Then $C$ is periodic with vertex
$v = c_1 + c_2 \in \Hi$.  By Lemma~\ref{Sebastiens lemma}, $S$ is
also periodic with vertex $v$.  Choose any hyperplane $\Sigma'$
of $\Sigma$ not containing $v$ and put $S' = S \cap \Sigma'$.
Then $S = \phi(S')$ (with respect to the vertex $v$). Since $S$ is
complete in $\Sigma$, it follows that $S'$ is complete in
$\Sigma'$.  Also the hyperplane $K_C \cap \Sigma'$ of $\Sigma'$
is a tangent hyperplane for $S'$ and thus $|S'|=2^{n-2}+1$. Hence
$|S| = 2|S'| = 2(2^{n-2}+1)=2^{n-1}+2$.
\end{proof}

\section{$|C|=3$}\label{triple threat}
  Now suppose that $S$ is a complete cap meeting a hyperplane
$K_C$ in 3 points.  Choose a codimension 2 subspace $\Hi \cong
\PG(n-2,2)$ contained in $K_C$ and disjoint from $S$. Write $C =
\{c_1,c_2,c_3\}$. Then $F = \{c_1+c_2,c_1+c_3,c_2+c_3\}$ is a
line in $\Hi$. Put $c_0 = c_1+c_2+c_3$. Then $\widetilde{F} = F
\cup C \cup \{c_0\} \cong \PG(2,2)$ and $\widehat{C} := C \cup
\{c_0\} \cong \AG(2,2)$. Thus $|C|=3$ implies that we are in the
case described in Section~\ref{family} with $r=2$.

  Thus by Theorem~\ref{sufficient conditions} (and
Remark~\ref{t or u odd}) we have maximal caps with $|C|=3$ of all
cardinalities $3 + 2m + 4(2^{n-3}-m) = 2^{n-1} - 2m + 3$ where $1
\leq m \leq 2^{n-3}-3$ or $m=2^{n-3}-1$, i.e., of cardinalities
$2^{n-2}+5, 2^{n-2}+9, 2^{n-2}+11, 2^{n-2}+13,\ldots, 2^{n-1}+1$
for $n \geq 4$.

Now we show that there for all $n \geq 3$ there exist complete
caps with $|C|=3$ and with $|A(i)|=|B(i)|=1$ for all
$i=1,2,\ldots,2^{n-3}$.  By Proposition~\ref{partition criterion}
we must find a partition of $\AG(n-2,2)$ into 4 subsets indexed by
$\AG(2,2)$ satisfying Equation~\ref{partition condition}.  There
are very many ways to do this.  We give one construction here.  We
proceed by induction on $k=n-2$. Since we have assumed that $r
\geq 2$ the first case is $k=4$.  We exhibited such a partition
for $k=4$ in Example~\ref{partition example}. Now we will
inductively construct from this example a partition for all
values of $k \geq 5$.  Suppose we have a partition of $\AG(k,2)$
into four sets $\b{A_{\emptyset}}$, $\b{A_0}$, $\b{A_1}$ and
$\b{A_{01}}$ where we further suppose that
$\b{A_{\emptyset}}=\emptyset$. Write $\PG(k+1,2) = \phi(\PG(k,2))$
with respect to some vertex $z \in PG(k+1,2) \setminus \PG(k,2)$.
Then $\AG(k+1,2) = \phi(\AG(k,2))$. Let $\wt{H} = \PG(k+1,2)
\setminus \AG(k+1,2) = \phi(H)$.  We partition $\AG(k+1,2)$ into
4 sets: $\AG(k+1,2) = \wt{A_\emptyset} \sqcup \wt{A_0} \sqcup
\wt{A_1} \sqcup \wt{A_{01}}$ by defining $\wt{A_w} :=
\phi(\b{A_w}) = \b{A_w} \sqcup (z+\b{A_w})$.
 This partition almost satisfies Equation~\ref{partition condition}.
However the point $z \in \wt{H}$ is not contained in any of the
three sets $\wt{A_0} + \wt{A_1}$, $\wt{A_0} + \wt{A_{01}}$ and
$\wt{A_1} + \wt{A_{01}}$. To overcome this defect we modify the
partition slightly. Note that every point of $\wt{H}$ different
from $z$ lies on at least two secants of each type $u \in
\AG(2,2)$.  Consider the three points of $\AG(k,2)$: $a_0 := 6$,
$a_1 := 236$ and $a_{01} := 26$ and the line in $H$: $z_{01} :=
a_0 + a_1$, $z_1 := a_0 + a_{01}$ and $z_0 := a_1 + a_{01}$. Each
of the three points of this line lies on secants of all types
which do not utilize $a_0, a_1$ nor $a_{01}$. Specifically we have
 $z_{01} = 23 = 256 + 356 = 2346 + 46 = 456 + 23456$,
 $z_1 = 2 = 2456 + 456 = 246 + 46 = 356 + 56$ and
 $z_0 = 3 = 256 + 2356 = 346 + 46 = 356 + 56$.
We shift $z+a_0$ into $\wt{A_1}$, $z+a_1$ into $\wt{A_{01}}$ and
$z+a_{01}$ into $\wt{A_0}$.  This modified partition now
satisfies Equation~\ref{partition condition}.  Also note that we
have not shifted any of the points $a_0, a_1, a_{01}$ nor indeed
any of the points of $\AG(k,2)$ and so our induction proceeds
using the same points $a_0, a_1, a_{01}$ at each stage.

  In conclusion we see that there are complete caps, $S$ in $\PG(n,2)$
with $|C|=3$ of cardinality $2^{n-2}+3$ for all $n \geq 7$. Thus
we have the following.

\begin{prop}
  Suppose $n \geq 4$ and that $S$ is a complete cap in $\PG(n,2)$.
Further suppose that $S$ meets a hyperplane $K_C$ of $\PG(n,2)$ in
exactly 3 points. Then $|S|$ is one of the numbers: $2^{n-2}+5,
2^{n-2}+9, 2^{n-2}+11, 2^{n-2}+13,\ldots, 2^{n-1}+1$ or, if $n
\geq 7$, $|S|$ may also be $2^{n-2}+3$ or $2^{n-2}+7$.  Moreover
complete caps of all these sizes exist for all corresponding
values of $n$.
\end{prop}

\section{$|C|=4$}\label{Four Slice}
  In this section we assume that $S$ is a complete cap meeting
a hyperplane $K_C$ in 4 points.  As usual we choose a codimension
2 subspace $\Hi \cong \PG(n-2,2)$ contained in $K_C$ and disjoint
from $S$.  Note that when $|K_C \cap S|=4$ there are two
possibilities; either these four points lie in a Fano plane or
they do not. Now if they do lie in a Fano plane they form a
periodic set and it follows as in Section~\ref{double plane} that
if $S$ complete then $|S| = 2^{n-1}+4$ and $S$ is large.

  It remains to consider the case where $S \cap K_C$ consists of four
points spanning a three dimensional projective subspace of
$K_C$.  In order to determine all solutions to
Equation~\ref{coset equations} we first assume that this three
dimensional subspace is $K_C$ and thus that $n=4$. Without loss
of generality, $C = \{04,14,24,34\}$, $K_C = (01234)^\perp$, $K_B
= (4)^\perp$ and $K_A = (0123)^\perp$.  Thus $\Hi =
\{01,02,03,12,13,23,\allowbreak 0123\}$ and $H_A = \HA{1} =
\{4,014,024,034,124,134,234,01234\}$.

  Let $G$ denote the group of projectivities of $\Sigma=\PG(4,2)$
which stabilize, $C$, $K_A$, $K_B$ and $K_C$.  Then $G$ contains
the subgroup $L$ of projectivities which permute the points
$0,1,2$ and $3$ and which also fix the point $4$.  The group $G$
also contains the projectivities $f$ and $h$ defined by: $f(0)=0$,
$f(1)=0$, $f(2)=0124$, $f(3)=0134$ and $f(4)=014$ and $h(0)=123$,
$h(1)=023$, $h(2)=013$, $h(3)=012$, and $h(4)=01234$.

  Suppose that $A$ and $B$ satisfy Equation~\ref{coset equations}
and further that $A$ and $B$ do not correspond to either trivial
or singleton solutions of that equation. Thus $|A| \geq 2$ and
$|B|\geq 2$.  Since $G$ acts transitively on $H_A$ we may assume
without loss of generality that $01234 \in A$.  Now the subgroup
of $G$ fixing $01234$ contains the group $L$.  The action of $L$
decomposes the remaining points of $H_A$ into two orbits: $\{4\}$
and $\{014,024,034,124,134,234\}$. Thus without loss of
generality we have two cases to consider: (i) $A \supset
\{01234,4\}$ and (ii) $A \supset \{01234,014\}$. These two cases
are distinguished by the fact that in the former case the point
$z_0 := c_1 + c_2 + c_3 + c_3 = 04 + 14 + 24 + 34 = 0123$ lies on
the secant through the two given points of $A$ while in the latter
case it does not.

  In the first case, we see that
$B'=A+C \supset\{01234,4\} + C = H_B$ and thus we have the
trivial solution $A=H_A$ and $B=\emptyset$.  In the second case
$B'=A+C \supset\{01234,014\} + C = H_B \setminus \{2,3\}$. Since
$|B| \geq 2$, we must have $B=\{2,3\}$ and $A=\{01234,014\}$.

  Thus we find only one new solution to Equation~\ref{coset equations}.
This solution is determined by $A(i) = \{a_1,a_2\}$, $a_1 + a_2
\neq z_0$, and $B(i) = \{b_1,b_2\} = \HB{i} \setminus (A(i)+C)$
where $z_0 = c_1 + c_2 + c_3 + c_4$.  Note that $b_1 + b_2 = a_1
+ a_2$.  Thus we have five solutions in total: the two trivial
solutions, the two singleton solutions and this new solution.
Notice that for the new solution and for both of the singleton
solutions we have $z_0 \notin A(i) \oplus A(i)$ and $z_0 \notin
B(i) \oplus B(i)$.  Since $z_0 \notin C \oplus C$, if $S$ is
complete there must be at least one $i$ for which the trivial
solution occurs.  In particular, we may apply
Theorem~\ref{sufficient conditions} to obtain complete caps $S$
with $|C|=4$ and $C$ not periodic.

  Suppose that $S$ is a complete cap with $|C|=4$.  If $C$ is
periodic then $|S|=2^{n-1}+4$.  Suppose then that $C$ is not
periodic and that $m$ of the pairs of cosets $\HA{i},\HB{i}$
utilize the new solution, $s$ of them utilize singleton solutions
and $t+u = 2^{n-4}-s-m$ of them utilize the trivial solution. Note
that if $s=0$ then we must have $m \geq 2$ in order that
$\widehat{C} \setminus C \subseteq \cup_{i=1}^{2^{n-r-1}}
\left(A(i) + B(i)\right)$. The cardinality of $S$ is then given
by $|S| = 4 + 4m + 5s + 8(2^{n-4}-s-m) = 2^{n-1} -3(s+m) -m +4$
where $0 \leq m,s \leq 2^{n-4}-1$ and $1 \leq m+s \leq 2^{n-4}-1$
(and $(m,s) \neq (1,0)$).  Thus if $|C|=4$ and $C$ is not periodic
then $n \geq 5$ and we find complete caps $S$ of all cardinalities
satisfying $2^{n-2}+ 8 \leq |S| \leq 2^{n-1}-2$.  In summary we
have the following.

\begin{prop}
  Suppose $n \geq 5$ and that $S$ is a complete cap in $\PG(n,2)$.
Further suppose that $S$ meets a hyperplane $K_C$ of $\PG(n,2)$ in
exactly 4 points. Then $2^{n-2}+ 8 \leq |S| \leq 2^{n-1}-2$ or
else $|S|=2^{n-1}+1$. Moreover except for $n=6$ for sizes 29 and
30, complete caps of all these sizes and this structure exist for
all $n \geq 5$.
\end{prop}

\begin{acknowledgements}
I thank Aiden Bruen for many useful conversations.  This research
was supported by grants from ARP and NSERC.
\end{acknowledgements}

 \medskip
\begin{flushleft}
\begin{tabular}{ll}
 \begin{tabular}{l}
         David L. Wehlau\\
         Department of Mathematics and Computer Science\\
         Royal Military College\\
         Kingston, Ontario\\
         Canada  K7K 7B4\\
         {\tt wehlau@rmc.ca}
 \end{tabular}
\end{tabular}
 \end{flushleft}

\begin{thebibliography}{99}

\bibitem{BBW} A.E. Brouwer, A. A. Bruen, and D.L. Wehlau,
{There exist caps which block all subspaces of fixed codimension
in $\PG(n,2)$}, \emph{J.C.T. A} \textbf{73}, \#1 (1996) 168--169.

\bibitem{BHW} Aiden A. Bruen, Lucien Haddad and David L. Wehlau,
{Binary Codes and Caps}, \emph{J. Combin. Des.} \textbf{6}, \#4
(1998) 275--284.


\bibitem{BW} Aiden A. Bruen and David L. Wehlau,
{Long Binary Linear Codes and Large Caps in Projective Space},
\emph{D.C.C.} \textbf{17}, (1999) 37--60.

\bibitem{Black/White lift} Aiden A. Bruen and David L. Wehlau,
{New Codes From Old; A New  Geometric Construction} (to appear in
\emph{J.C.T. A}).

\bibitem{DT} A.A. Davydov and L.M. Tombak,
{Quasiperfect Linear Binary Codes with Distance 4 and
Complete Caps in Projective Geometry}, \emph{Problems of Information
Transmission}~\textbf{25} No.~4 (1990) 265--275.

 \end{thebibliography}
 \end{document}